**A new method for optimal control of Volterra integral equations**

S. A. Belbas
Mathematics Department
University of Alabama
Tuscaloosa, AL. 35487-0350. USA.

e-mail: SBELBAS@GP.AS.UA.EDU

Abstract. We formulate and analyze a new method for solving optimal control problems for systems governed by Volterra integral equations. Our method utilizes discretization of the original Volterra controlled system and a novel type of dynamic programming jn which the Hamilton-Jacobi function is parametrized by the control function (rather than the state, as in the case of ordinary dynamic programming). We also derive estimates for the computational cost of our method.

## 1. Introduction.

The classical theory of optimal control was originally developed to deal with systems of controlled ordinary differential equations. It has been understood that many physical, technological, biological, and socio-economic problems cannot be adequately described by ordinary differential equations, and other mathematical models, including systems with memory, distributed systems, and other types of systems, have been added to the arsenal of the theory of optimal control. A broad category of systems can be described by Volterra integral equations.

The simplest form of a controlled Volterra integral equation is

$$x(t) = x_0(t) + \int_0^t f(t,s,x(s),u(s))\,ds \qquad \text{--- (1.1)}$$

In this system, x(t) is the n-dimensional state function, and u(t) is the m-dimensional control function. For the purposes of this exposition, we postulate that f be continuous with respect to all variables and uniformly Lipschitz with respect to x. For the purposes of describing the necessary conditions that are briefly reviewed in this section, the admissible control functions are continuous functions with values in a compact set **U**, $\mathbf{U} \subseteq \mathbb{R}^n$ . In certain parts of the overall theory of optimal control for Volterra integral equations, the class of admissible control functions can be more general, for example it may consist of bounded measurable or p-integrable functions. On the other hand, as it will be explained below, for some of the results of the present paper it is necessary to further restrict the class of admissible controls and postulate Lipschitz continuity.

Volterra integral equations arise in a wide variety of applications. In fact, it seems that, with the exception of the simplest physical problems, practically every situation that can be modelled by ordinary diffrential equations can be extended to a model with Volterra integral equations. For example, a general ODE system of interacting biological populations, of the form

$$\frac{dx_i(t)}{dt} = \sum_{j,k} a_{ijk} x_j(t) x_k(t) + \sum_j b_{ij} x_j(t);\ x_i(0) = x_{i,0}$$

can be extended to an integro-differential system

$$\frac{dx_i(t)}{dt} = \sum_{j,k} a_{ijk} x_j(t) x_k(t) + \sum_j b_{ij} x_j(t) +$$

$$+ \int_0^t \left\{ \sum_{j,k} A_{ijk}(t,s) x_j(s) x_k(s) + \sum_j B_{ij}(t,s) x_j(s) \right\} ds;$$

$$x_i(0) = x_{i,0}$$

Indeed, some related extensions have already been considered in [CU] and in other works. In turn, every integrodifferential system of the form

$$\frac{dx(t)}{dt} = g\left( t, x(t), \int_0^t f(t,s,x(s))\, ds \right); \; x(0) = x_0$$

can be reduced to a system of Volterra integral equations by setting $y(t) = \int_0^t f(t,s,x(s)\, ds$, then the integrodifferential equation becomes $x(t) = x_0 + \int_0^t g(s,x(s),y(s))\, ds$, so that we get a system of Volterra integral equations in the unknowns (x(t), y(t)).

Problems in mathematical economics also lead to Volterra integral equations. The relationships among different quantities, for example between capital and investment, include memory effects (e.g. the present stock of capital depends on the history of investment strategies over a period of time, cf. [KM]), and the simplest way to describe such memory effects is through Volterra integral operators.

Now we return to the general model of state dynamics (1.1).

An optimal control problem for (1.1) concerns the minimization of a cost functional

$$J := F_0(x(T)) + \int_0^T F(t,x(t),u(t))\, dt \qquad \text{--- (1.2)}$$

The theory of optimal control of ordinary differential equations has two main methods: the extremum (usually called "maximum") principle of Pontryagin and his coworkers, and the method of dynamic programming. The metod of dynamic programming is particularly useful as it provides sufficient conditions for optimality. However, the nature of controlled Volterra equations is not, at first glance, conducive to the application of dynamic programming methods. If the state x(t) is known at some particular time t, and a control function is specified over an interval $(t, t+\delta t]$, these two bits of information are not enough for the determination of the solution of (1.1) over the interval $(t, t+\delta t]$. By contrast, for ordinary differential equations, it is always true that, given x(t) and a control function over $(t, t+\delta t]$, the trajectory over $(t, t+\delta t]$ can be determined by solving an initial value problem for an ordinary differential equation with initial time t. For these reasons, optimal control problems for Volterra integral equations have been traditionally treated by extensions of Pontryagin's extremum principle. The related results are found in a number of papers, including [M, S, V, NW]; an approach based on direct variational methods, but still utilizing necessary conditions for optimality, may be found in [B].

The co-state $\psi(t)$ for the problem consisting of (1.1) and (1.2) satisfies the following adjoint equation, which is the counterpart of Hamiltonian equations (see, e.g., [S, V]):

$$\psi(t) = \nabla_x F(t, x(t), u(t)) + \left(\nabla_x F_0(x(T))\right)\nabla_x f(T, t, x(t), u(t)) + \\ + \int_t^T \psi(s)\nabla_x f(s, t, x(t), u(t))\,ds$$

--- (1.3)

The state x(.) is an n-dimensional column vector; the function f takes values that are n-dimensional column vectors; the co-state $\psi(.)$ is an n-dimensional row vector. The gradient, with respect to x, of a scalar-valued function is an n-dimensional row vector; $\nabla_x f$ is a matrix with elements $\frac{\partial f_i}{\partial x_j}$ (each $\frac{\partial f_i}{\partial x_j}$ is the element in the $i^{th}$ row and $j^{th}$ column of $\nabla_x f$ ). It should be noted that the adjoint equations in [S] and [V] are not exactly identical to (1.3), due to the fact that these authors do not use exactly the same cost functional as (1.2); however, these differences do not require substantially different proofs, and for that reason we present (1.3) without proof.

The Hamiltonian is defined by

$$H(t, x, \psi(.), u) := F(t, x, u) + \left(\nabla_x F_0(x(T))\right) f(T, t, x, u) + \int_t^T \psi(s) f(s, t, x, u)\,ds$$

--- (1.4)

The extremum principle takes the following form [V]: under certain smoothness conditions, an optimal control $u^*(.)$ and the corresponding trajectory and co-state, $x^*(.)$ and $\psi^*(.)$, respectively, satisfy, for almost all t,

$$H(t,x^*(t),\psi^*(.),u^*(t)) \le H(t,x^*(t),\psi^*(.),u) \quad \forall u \in \mathbf{U}$$
--- (1.5)

In this paper, we shall need the concept of what we shall term the <u>relevant values</u> of the state x(.). This is not standard terminolgy, but is is useful for our purposes.
Under certain conditions, it is possible to estimate the range of the solution x(t), without relying on actually solving (1.1). For example, if the function f has linear growth rate with respect to x, i.e. if

$$|f(t,s,x,u)| \le G_1|x| + G_0, \text{ for all } t,s: 0 \le s \le t \le T, \text{ for all } u \in \mathbf{U}$$
--- (1.6)

At first, we do not specify the set of values of x for which (1.6) should hold; we tentatively carry out the calculations as if (1.6) were true for all x in $\mathbb{R}^n$, until we find an estimate for the set $X_{rel}$ of relevant values of x(.), then we go back and postulate that (1.6) should hold for all x in $X_{rel}$.

Then it follows from (1.1) that every solution x(.) satisfies

$$|x(t)| \le \|x_0\|_\infty + TG_0 + G_1\int_0^t |x(s)|ds; \ \|x_0\|_\infty := \sup_{0\le t\le T} |x_0(t)|$$
--- (1.7)

and therefore, by Gronwall's inequality,

$$|x(t)| \le (\|x_0\|_\infty + TG_0)\exp(G_1T) \quad \forall \in [0,T]$$
--- (1.8)

In this case, the set $X_{rel}$ of relevant values of x(.) is

$$X_{rel} := \{x \in \mathbb{R}^n : |x| \le (\|x_0\|_\infty + TG_0)\exp(G_1T)\}$$
--- (1.9)

Thus the conclusion is that, if (1.6) holds for all x in the set $X_{rel}$ given by (1.9), then the solution x(t) always lies in that $X_{rel}$.

It should be emphasized that (1.9) is an example of one case of a set of relevant values, and not the definition of $X_{rel}$ in general.

Another example of estimating $X_{rel}$ is the case in which f satisfies a Lipschitz condition,

$$|f(t,s,x_1,u)-f(t,s,x_2,u)|\le L_{f,1}|x_1-x_2|,\ \text{for all } t,s: 0\le s\le t\le T, \text{ for all } u\in \mathbf{U} \quad \text{--- (1.10)}$$

and for a particular value of x, which we may take, without loss of generality, to be 0, the function f remains bounded:

$$|f(t,s,0,u)|\le G_2\ ,\ \text{for all } t,s: 0\le s\le t\le T, \text{ for all } u\in \mathbf{U} \quad \text{--- (1.11)}$$

As before, we work at first as if (1.10) were true for all $x_1$, $x_2$ in $\mathbb{R}^n$, then we determine a suitable set $X_{rel}$, and then we return and postulate that (1.10) should hold for all $x_1$, $x_2$ in $X_{rel}$.
Now, we have

$$|x(t)|-\int_0^t |f(t,s,0,u(s))|ds \le |x(t)-\int_0^t f(t,s,0,u(s))ds| \le$$
$$\le |x_0(t)|+L_{f,1}\int_0^t |x(s)|ds \quad \text{--- (1.12)}$$

thus

$$|x(t)|\le G_2T+\|x_0\|_\infty+L_{f,1}\int_0^t |x(s)|ds \quad \text{--- (1.13)}$$

and, again by Gronwall's inequality,

$$|x(t)| \le (G_2 T + \| x_0 \|_\infty) \exp(L_{f,1} T) \ \forall t \in [0,T] \qquad \text{--- (1.14)}$$

so that, in this case, we can take

$$X_{rel} := \{x \in \mathbb{R}^n : |x| \le (G_2 T + \| x_0 \|_\infty) \exp(L_{f,1} T)\} \qquad \text{--- (1.15)}$$

Thus, if (1.10) and (1.11) hold for all $x_1$, $x_2$ in the set $X_{rel}$ given by (1.15), then all values of x(t) lie in that set $X_{rel}$.

It should be noted that there are many other possibilities of finding examples of $X_{rel}$ under suitable assumptions, but, in this paper, we are not interested in exhausting this topic.

In the rest of this paper, we shall assume the existence of a bounded set $X_{rel} \subseteq \mathbb{R}^n$, without specifying how that set has been determined.

## 2. The discrete Volterra control problem.

Our approach will be to aproximate the original Volterra control problem by a sequence of analogous control problems for discrete Volterra equations. For this reason we need to have a method for solving optimal control problems for discrete Volterra equations. Thus we consider, in this section, the controlled Volterra equation in discrete time:

$$x(i) = x_0(i) + \sum_{j=0}^{i-1} \varphi(i, j, x(j), u(j)), \; 1 \le i \le N \qquad \text{--- (2.1)}$$

The discrete optimal control problem concerns the minimization of a functional J given by

$$J := \sum_{i=0}^{N-1} \Phi(i, x(i), u(i)) + \Phi_0(x(N)) \qquad \text{--- (2.2)}$$

In order to apply a suitable variant of the dynamic programing method to the problem consisting of (2.1) and (2.2), we need to build a parametrization of this optimal control problem. The expression "suitable variant" refers to the fact that, for the problem under consideration, the value function needs to be parametrized by current time and history of the control up to the current time, whereas in classical dynamic programming the value function is parametrized by current time and current value of the state. The memory effect of Volterra equations necessitates this seemingly unorthodox parametrization.
We set

$$\tilde{u}_{(i)} := (u(j) : 0 \le j \le i-1) \qquad \text{--- (2.3)}$$

If $x(k; i, \tilde{u}_{(i)})$, $0 \le k \le i$ the solution of the discrete Volterra equation

$$x(k) = x_0(k) + \sum_{j=0}^{k-1} \varphi(k, j, x(j), u(j)), \; 0 \le k \le i-1 \qquad \text{--- (2.4)}$$

We define

$$x_{[i]} \equiv x_{[i]}(i, \tilde{u}_{(i)}) := (x(k; i, \tilde{u}_{(i)}) : 0 \le k \le i) \qquad \text{--- (2.5)}$$

The <u>concatenation</u> of $\beta = (\beta(j) : 0 \le j \le i-1)$ and $\gamma = (\gamma(j) : i \le j \le N-1)$ is defined as

$$\lambda \equiv \beta \otimes \gamma;\ \lambda(j) = \beta(j) \text{ for } j = 0,1,\ldots,i-1; \lambda(j) = \gamma(j) \text{ for } j = i, i+1, \ldots, N-1 \qquad \text{--- (2.6)}$$

The restriction of a control to indices that exceed $i-1$ will be denoted by $u_{<i>}$:

$$u_{<i>} = (u(j) : i \le j \le N-1) \qquad \text{--- (2.7)}$$

The cost functional J is parametrized as

$$J_{i,\beta}(u_{<i>}) := \sum_{j=i}^{N-1} \Phi(j; x(j; i, \beta \otimes u_{<i>})), u(j)) + \Phi_0(x(N; i, \beta \otimes u_{<i>})) \qquad \text{--- (2.8)}$$

where $\beta = (\beta(j) : 0 \le j \le i-1)$, $u_{<i>} = (u(j) : i \le j \le N-1)$, and $x(k; i, \beta \otimes u_{<i>}),\ i \le k \le N-1$ solves

$$x(k) = x_0(k) + \sum_{j=0}^{i} \varphi(k, j, x(j; i, \beta), \beta(j)) + \sum_{j=i+1}^{k-1} \varphi(k, j, x(j), u(j)) \qquad \text{--- (2.9)}$$

In particular, we note that $x(i; i, \beta \otimes u_{<i>}) = x(i; i, \beta)$.

The <u>value function</u> $V(i, \beta)$, $\beta = (\beta(j) : 0 \le j \le i-1)$, is defined in terms of the parametrization (2.8):

$$V(i, \beta) := \min_{u_{<i>}} J_{i,\beta}(u_{<i>}) \qquad \text{--- (2.10)}$$

For $i=0$, the collection $(\beta(j):0\le j\le i-1)$ is empty; we therefore use a symbolic "empty set" $\varnothing$ in the function V, and that function becomes, when $i=0$, $V(0,\varnothing)$.

For $\xi\in U$, we identify $\beta\otimes\xi$ with the control $\hat{\beta}:=(\beta(0),\beta(1),...,\beta(i-1),\xi)$ so that $\hat{\beta}(i)=\xi$. With this notational convention, the dynamic programming equations for V are

$$V(i,\beta)=\min_{\xi\in\mathbf{U}}\{V(i+1,\beta\otimes\xi)+\Phi(i,x(i;i,\beta),\xi)\};\ V(N,\beta)=\Phi_0(x(N;N,\beta)) \quad \text{--- (2.11)}$$

Next, we prove the necessity and sufficiency of (2.10).

We have:

Theorem 2.1. Eq. (2.11) is necessary for optimality, i.e. if $V(i,\beta)$ is defined by (2.10), then it satisfies (2.11).

Proof: According to (2.10), we have, for every control $u_{<i+1>}$,

$$V(i+1,\beta\otimes\xi)\le J_{i+1,\beta\otimes\xi}(u_{<i+1>}) \quad \text{--- (2.12)}$$

Therefore, for every $\xi\in\mathbf{U}$, we have

$$\begin{aligned}&\min_{\xi\in\mathbf{U}}\{V(i+1,\beta\otimes\xi)+\Phi(i,x(i;i,\beta),\xi)\}\le V(i+1,\beta\otimes\xi)+\Phi(i,x(i;i,\beta),\xi)\le\\&\le J_{i+1,\beta\otimes\xi}(u_{<i+1>})+\Phi(i,x(i;i,\beta),\xi)=J_{i,\beta}(\xi\otimes u_{<i+1>})\end{aligned} \quad \text{--- (2.13)}$$

Since every $u_{<i>}$ can be represented as $\xi\otimes u_{<i+1>}$ for some $\xi\in\mathbf{U}$, we have

$$\begin{aligned}V(i,\beta)&=\min_{u_{<i>}}J_{i,\beta}(u_{<i>})=\min_{\xi\in\mathbf{U}}\ \min_{u_{<i+1>}}\ J_{i,\beta}(\xi\otimes u_{<i+1>})=\\&=\min_{\xi\in\mathbf{U}}\{V(i+1,\beta\otimes\xi)+\Phi(i,x(i;i,\beta),\xi)\}\end{aligned} \quad \text{--- (2.14)}$$

///.

<u>Theorem 2.2.</u> Eq. (2.11) is sufficient for optimality, i.e. the solution of (2.11) satisfies (2.10). If a control function $u^*(.)$ satisfies $V(i,\tilde{u}^*_{(i)}) = J_{i,\tilde{u}^*_{(i)}}(u^*_{<i>})$ for all $i \in \{0,1,2,...,N\}$, then $u^*(.)$ is an optimal control, i.e. $J(u^*(.)) \le J(u(.))$ for every admissible control function u(.).

<u>Proof:</u> The proof that (2.11) implies (2.10) uses backward induction. If the statement is true for i+1, we shall show that it must be true for i. By the induction hypothesis, we have

$$V(i+1,\beta \otimes \xi) = \min_{u_{<i+1>}} J_{i+1,\beta\otimes\xi}(u_{<i+1>}) \qquad \text{--- (2.15)}$$

As in the proof of theorem 2.1, every $u_{<i>}$ can be represented as $u_{<i>} = \xi \otimes u_{<i+1>}$. Then we have

$$J_{i,\beta}(\xi \otimes u_{<i+1>}) = J_{i+1,\beta\otimes\xi}(u_{<i+1>}) + \Phi(i, x(i;i,\beta),\xi) \qquad \text{--- (2.16)}$$

Consequently

$$\begin{aligned}
&\min_{\xi\otimes u_{<i+1>}} J_{i,\beta}(\xi \otimes u_{<i+1>}) = \min_{\xi\otimes u_{<i+1>}} \{J_{i+1,\beta\otimes\xi}(u_{<i+1>}) + \Phi(i, x(i;i,\beta),\xi)\} = \\
&= \min_{\xi} \min_{u_{<i+1>}} \{J_{i+1,\beta\otimes\xi}(u_{<i+1>}) + \Phi(i, x(i;i,\beta),\xi)\} = \\
&= \min_{\xi} \{\{\min_{u_{<i+1>}} J_{i+1,\beta\otimes\xi}(u_{<i+1>})\} + \Phi(i, x(i;i,\beta),\xi)\} = \min_{\xi}\{V(i+1,\beta\otimes\xi) + \Phi(i, x(i;i,\beta),\xi)\} = \\
&= V(i,\beta)
\end{aligned} \qquad \text{--- (2.17)}$$

thus

$$V(i,\beta) = \min_{u_{<i>}} J_{i,\beta}(u_{<i>}) \qquad \text{--- (2.18)}$$

For i=N, the statement is true since $J_{N,\beta}(u_{<N>}) \equiv J_{N,\beta}(\varnothing) = \Phi_0(x(N;N,\beta)) = V(N,\beta)$. Thus the backward induction is complete.

For the second asertion of this theorem, suppose $u^*(.)$ satisfies $V(i,\tilde{u}^*_{(i)}) = J_{i,\tilde{u}^*_{(i)}}(u^*_{<i>})$ for all $i \in \{0,1,2,...,N\}$. Then, in particular, for i=0, we have $V(0,\varnothing) = J_{0,\varnothing}(u^*_{<0>})$. It is a simple consequence of our notational convention that, for every admissible control function u(.), we have $J_{0,\varnothing}(u_{<0>}) \equiv J(u)$. Therefore, $V(0,\varnothing) = J(u^*(.))$. At the same time, since V satisfies (2.9), we have $V(0,\varnothing) \le J_{0,\varnothing}(u_{<0>}) \equiv J(u(.))$ for every admissible control function u(.). Therefore $J(u^*(.)) \le J(u(.))$ for every admissible control function u(.). ///.

The next question is how to find a control function $u^*(.)$ that satisfies $V(i,\tilde{u}^*_{(i)}) = J_{i,\tilde{u}^*_{(i)}}(u^*_{<i>})$ for all $i \in \{0,1,2,...,N\}$. This is done by forward recursion. For i=0, $u^*(0)$ is found as a value of ξ that is a minimizer of $V(1,\xi) + \Phi(0,x_0(0),\xi)$, or, equivalently, as a solution of $V(0,\varnothing) = V(1,u^*(0)) + \Phi(0,x_0(0),u^*(0))$. Inductively, if $u^*_{(i)}$ has been determined, then $u^*(i)$ is found as a value of ξ that minimizes $V(i+1,u^*_{(i)} \otimes \xi) + \Phi(i,x(i;i,u^*_{(i)}),\xi)$, or, equivalently, as a solution of $V(i,u^*_{(i)}) = V(i+1,u^*_{(i)} \otimes u^*(i)) + \Phi(i,x(i;i,u^*_{(i)}),u^*(i))$. We have the following:

<u>Theorem 2.3.</u> If a control function $u^*(.)$ is constructed so as to satisfy

$$V(0,\varnothing) = V(1,u^*(0)) + \Phi(0,x_0(0),u^*(0));$$
$$V(i,u^*_{(i)}) = V(i+1,u^*_{(i)} \otimes u^*(i)) + \Phi(i,x(i;i,u^*_{(i)}),u^*(i)) \text{ for } 1 \le i \le N-1 \quad \text{--- (2.19)}$$

then it also satisfies

$$V(i,\tilde{u}^*_{(i)}) = J_{i,\tilde{u}^*_{(i)}}(u^*_{<i>}) \text{ for all } i \in \{0,1,2,...,N\} \quad \text{--- (2.20)}$$

<u>Proof:</u> We denote by $x^*(.)$ the solution of (2.1) that corresponds to the control function $u^*(.)$. We use (2.11) and the fact that $u^*_{<N>} \equiv \varnothing$ to obtain

$$V(N,u^*_{(N)}) = \Phi_0(x^*(N)) \equiv J_{N,u^*_{(N)}}(\varnothing) \equiv J_{N,u^*_{(N)}}(u^*_{<N>})$$

--- (2.21)

thus the wanted assertion is true for $i = N$. We use backward induction: assuming that the wanted assertion is true for i+1, we shall show that it must be true for i. By the induction hypothesis, we have

$$V(i+1,u^*_{(i+1)}) = J_{i+1.u^*_{(i+1)}}(u^*_{<i+1>})$$

--- (2.22)

Since $V(i,u^*_{(i)}) = V(i+1,u^*_{(i)} \otimes u^*(i)) + \Phi(i,x(i;i,u^*_{(i)}),u^*(i))$ and $u^*_{<i>} = u^*(i) \otimes u^*_{<i+1>}$, we have

$$V(i,u^*_{(i)}) = J_{i+1.u^*_{(i+1)}}(u^*_{<i+1>}) + \Phi(i,x(i;i,u^*_{(i)}),u^*(i)) = J_{i,u^*_{(i)}}(u^*_{<i>})$$

--- (2.23)

thus the induction is complete. ///

Remark 2.1. The construction of an optimal control for our variant of dynamic programming for discrete Volterra equations differs in a substantial way from the "feedback" or "closed loop" controls that are obtained in ordinary dynamic programming (i.e. in dynamic programming for ordinary differential equations or finite-difference equations). In our case, each optimal value $u^*(i)$ depends on, among other things, the future optimal control policy $u^*_{<i>}$. This additional complication further contributes to Bellman's "curse of dimensionality". It is, of course, natural that, due to the memory effect of Volterra equations, the construction of optimal controls will be more complicated than in the case of ordinary differential equations or finite-difference equations. ///

Remark 2.2. The method of dynamic programming developed above can be modified to include the possibility of constraints of the type

$$u(i) \in \Xi(i,\tilde{u}_{(i)},\tilde{x}_{[i]})$$

--- (2.24)

where each set $\Xi(i,\tilde{u}_{(i)},\tilde{x}_{[i]})$ is a closed subset of **U**. In that case, the dynamic programming equations take the form

$$V(i,\beta) = \min_{\xi \in \Xi(i,\beta,\tilde{x}_{[i]}(i,\beta))} \{V(i+1, \beta \otimes \xi) + \Phi(i, x(i;i,\beta), \xi)\}; \; V(N,\beta) = \Phi_0(x(N;N,\beta)) \quad \text{--- (2.25)}$$

The proof that these dynamic programming equations are necessary and sufficient for optimality under the constraints $u(i) \in \Xi(i, \tilde{u}_{(i)}, \tilde{x}_{[i]})$ can be carried out as in the unconstrained case, and therefore we omit the details. ///

3. Results on discretization of controlled Volterra equations and cost functionals.

We consider an Euler discretization of the original Volterra controlled system. If $h \equiv \frac{T}{N}$ is the step size of the Euler discretization, we set $t_i := ih$, i=0, 1, 2, ..., N, $x_0^h(i) := x_0(t_i)$, $f^h(i,j,x,u) := f(t_i,t_j,x,u)$, $u^h(i) := u(t_i)$. The discretized controlled Volterra equation is

$$x^h(i) = x_0^h(i) + h\sum_{j=0}^{i-1} f^h(i,j,x^h(j),u^h(j)),\ 0 \le i \le N \qquad \text{--- (3.1)}$$

The cost functional J is discretized as

$$J^h(x(.),u(.)) := h\sum_{i=0}^{N-1} F(t_i,x(t_i),u^h(i)) + F_0(x(T)) \qquad \text{--- (3.2)}$$

We also consider the functional

$$J^h(x^h(.),u(.)) := h\sum_{i=0}^{N-1} F(t_i,x^h(i),u^h(i)) + F_0(x^h(N)) \qquad \text{--- (3.3)}$$

In (3.2), x(.) is the solution of the continuous Volterra integral equation (1.1), whereas in (3.3) $x^h(.)$ is the solution of the discretized Volterra equation (3.1). We have chosen the simplest numerical integration schemes in order to minimize the regularity assumptions that we need for the error estimates.

The existing literature on numerical solution of Volterra integral equations deals with simple (i.e. non-controlled) Volterra equations. The approximation of controlled Volterra equations involves additional ingredients, and for this reason we cannot simply invoke existing results, but instead we must prove all the results we need.

For the purpose of obtaining error estimates, it becomes necessary to restrict the class of admissible control functions to functions that satisfy a Lipschitz condition.

Definition 3.1. The set of admisible control functions $\mathbf{U}_{ad,Lip}(L)$ is defined as the set of all functions u(.) from [0, T] into **U** that satisfy a Lipschitz condition with fixed Lipschitz constant L: $|u(t_1)-u(t_2)|\leq L|t_1-t_2|$ for all $t_1,t_2$ in [0, T]. ///

We postulate:

(i). There is a compact subset $X_{rel}$ of $\mathbb{R}^n$ such that all values x(t) of solutions of the continuous problem (1.1) and also all values $x^h(t_i)$ of solutions of the discrete Volterra equation fall into the set $X_{rel}$. (In other words, $X_{rel}$ is a common set of relevant values for both the continuous and the discretized problem.

Also, we postulate the following properties for the functions f, F, and $F_0$, in addition to the previous conditions:

(ii). The function f is jointly Lipschitz in x and u, with Lipschitz constant $L_f$, uniformly in s and t:

$|f(t,s,x_1,u_1)-f(t,s,x_2,u_2)|\leq L_f\{|x_1-x_2|+|u_1-u_2|\}$ for all $x_1,x_2$ in $\mathbb{R}$ and all $u_1,u_2$ in **U**, for all (s, t) that satisfy $0\leq s\leq t\leq T$.

(iii). The functions f, $f_t$, $f_s$ are bounded:
$\max\{|f(t,s,x,u)|,|f_t(t,s,x,u)|,|f_s(t,s,x,u)|\}\leq C_0$, for all $x\in\mathbb{R}, u\in\mathbf{U}$, and (s, t) that satisfy $0\leq s\leq t\leq T$.

(iv). The function F is jointly Lipschitz in x and u, and has bounded derivative with respect to t:

$|F(t,x_1,u_1)-F(t,x_2,u_2)|\leq L_F\{|x_1-x_2|+|u_1-u_2|\}, \forall x_1,x_2$ in $X_{rel}$, $\forall u_1,u_2$ in **U**;
$|F_t(t,x,u)|\leq M_F, \forall x\in X_{rel}, \forall u\in\mathbf{U}$.

(v). The function $F_0$ is Lipschitz: $|F_0(x_1)-F_0(x_2)|\leq L_{F_0}|x_1-x_2|\ \forall x_1,x_2$ in $X_{rel}$.

Some explanations are in order about condition (i). When the set $X_{rel}$ is found either from (1.9) or from (1.15), under the appropriate conditions in each case, then the same set $X_{rel}$ contains also all values $x^h(t_i)$ of all solutions of the discretized Volterra equation. In the case of (1.9), under condition (1.6), we have, for the discretized problem,

$$|x^h(i)| \le \|x_0\|_\infty + h\sum_{j=0}^{i-1} (G_0 + G_1|x^h(j)|) = \|x_0\|_\infty + TG_0 + hG_1\sum_{j=0}^{i-1} |x^h(j)|$$

from which, by the discrete Gronwall inequality,

$$|x^h(i)| \le (\|x_0\|_\infty + TG_0)(1+hG_1)^{\frac{i}{h}} \le (\|x_0\|_\infty + TG_0)(1+hG_1)^{\frac{T}{h}} .$$

Since $(1+hG_1)^{\frac{T}{h}} \downarrow \exp(G_1T)$ as $h \to 0^+$, we conclude that every $x^h(i)$ is in the set $X_{rel}$ given by (1.9).
In the case of (1.15) under conditions (1.10) and (1.11), it can be proved, in a similar way, that every $x^h(i)$ is in the set $X_{rel}$ given by (1.15).
Consequently, condition (i), in its general form, is a reasonable condition that can be satisfied in specific cases.

We have:

Theorem 3.1. Under the conditions of section 2 and the conditions (i) and (ii) above, for every $u(.) \in \mathbf{U}_{ad,Lip}(L)$ , with $h = \frac{T}{N}$, the solution of (3.1) satisfies

$|x(ih) - x^h(i)| \le C_1 h$ for all $i = 1,2,...,N$, uniformly in N and $u(.) \in \mathbf{U}_{ad,Lip}(L)$

--- (3.4)

where $C_1$ is a constant that can be expressed in terms of L, $L_f$, and $C_0$.

Proof: We set

$$M_1 := \sup\{|f_s(t,s,x,u)| : 0 \le s \le t \le T, x \in X_{rel}, u \in \mathbf{U}\};$$
$$M_2 := \sup\{|f_t(t,s,x,u)| : 0 \le s \le t \le T, x \in X_{rel}, u \in \mathbf{U}\}$$

--- (3.5)

The conditions we have postulated lead to, among other things, a uniform bound on the time-derivative of x(.). We use a dot to denote the time-derivative of x. We have:

$$\dot{x}(t) = \dot{x}_0(t) + f(t,t,x(t),u(t)) + \int_0^t f_t(t,s,x(s),u(s))ds$$

thus

$$|\dot{x}(t)| \le |\dot{x}_0(t)| + |f(t,t,x(t),u(t))| + \int_0^t |f_t(t,s,x(s),u(s))|ds \le$$
$$\le M_0 + C_0 + M_1 t \le M_0 + C_0 + M_1 T \equiv M_x$$

--- (3.6)

We note the following fact from analysis: if φ is a differentiable function from [0, T] into $\mathbb{R}^d$, with components $\varphi_i, i = 1,2,...,d$, with derivative that satisfies $|\dot{\varphi}(t)| \le M_\varphi \ \forall t \in [0,T]$, where $|\cdot|$ denotes one of the $\ell^p$ norms on $\mathbb{R}^d$, $1 \le p \le \infty$, then for every two points $t_1, t_2$ in $\mathbb{R}^d$ we have $|\varphi(t_1) - \varphi(t_2)| \le \lambda_d M_\varphi |t_1 - t_2|$, where the coefficient $\lambda_d$ depends on the dimension d and the norm that we use. For the norm $|v|_p := \left( \sum_{i=1}^{d} |v_i|^p \right)^{\frac{1}{p}}$, $1 \le p < \infty$, we have $\lambda_d = d^{\frac{1}{p}}$, and for the norm $|v|_\infty := \max_{1\le i\le d} |v_i|$ we have $\lambda_d = 1$.

Next, we estimate the error of approximation $|x(ih) - x^h(i)|$ as follows:

$$|x(ih) - x^h(i)| \le \sum_{j=0}^{i-1} \int_{jh}^{(j+1)h} |f(ih,s,x(s),u(s)) - f(ih,jh,x^h(j),u(jh))|ds \le$$
$$\le \sum_{j=0}^{i-1} \int_{jh}^{(j+1)h} \{|f(ih,s,x(s),u(s)) - f(ih,jh,x(s),u(s))| +$$
$$+ |f(ih,jh,x(s),u(s)) - f(ih,jh,x(jh),u(jh))| +$$
$$+ |f(ih,jh,x(jh),u(jh)) - f(ih,jh,x^h(j),u(jh))|\} ds \le$$
$$\le \sum_{j=0}^{i-1} \int_{jh}^{(j+1)h} \{M_1 \lambda_n |s - jh| + L_f[|x(s) - x(jh)| + |u(s) - u(jh)|] +$$
$$+ L_f |x(jh) - x^h(j)|\} ds$$

--- (3.7)

The various terms on the right-hand side of the last inequality in (3.7) are estimated as follows:

$$\sum_{j=0}^{i-1}\int_{jh}^{(j+1)h} M_1\lambda_n\,|s-jh|\,ds \le \sum_{i=0}^{i-1} M_1\lambda_n h^2 = ih^2 M_1\lambda_n \le ThM_1\lambda_n\,;$$

$$\sum_{j=0}^{i-1}\int_{jh}^{(j+1)h} L_f\,|x(s)-x(jh)|\,ds \le \sum_{j=0}^{i-1}\int_{jh}^{(j+1)h} L_f\lambda_n M_x\,|s-jh|\,ds \le$$

$$\le \sum_{j=0}^{i-1}\int_{jh}^{(j+1)h} L_f\lambda_n M_x h\,ds = \sum_{i=0}^{i-1} L_f\lambda_n M_x h^2 = ih^2 L_f\lambda_n M_x \le ThL_f\lambda_n M_x\,;$$

$$\sum_{j=0}^{i-1}\int_{jh}^{(j+1)h} L_f\,|u(s)-u(jh)|\,ds \le \sum_{j=0}^{i-1}\int_{jh}^{(j+1)h} L_f M_2\lambda_m\,|s-jh|\,ds \le$$

$$\le \sum_{j=0}^{i-1}\int_{jh}^{(j+1)h} L_f M_2\lambda_m h\,ds = \sum_{i=0}^{i-1} L_f M_2\lambda_m h^2 = ih^2 L_f M_2\lambda_m \le ThL_f M_2\lambda_m$$

--- (3.8)

We set

$$\overline{M} := T[M_1\lambda_n + L_f\lambda_n M_x + L_f M_2\lambda_m]$$

--- (3.9)

Then (3.7) and (3.9) give

$$|\,x(ih) - x^h(i)\,| \le \overline{M}h + L_f h \sum_{j=0}^{i-1} |\,x(jh) - x^h(j)\,|$$

--- (3.10)

from which we obtain, via the discrete Gronwall inequality,

$$|\,x(ih) - x^h(i)\,| \le \overline{M}h(1+hL_f)^i \le \overline{M}h(1+hL_f)^N = \overline{M}h(1+hL_f)^{\frac{T}{h}} \le \overline{M}h\exp(L_f T)$$

--- (3.11)

which proves the assertion of the theorem, with $C_1 = \overline{M}\exp(L_f T)$. ///

Theorem 3.2. Under the above conditions, we have

$$| J(x(.),u(.)) - J^h(x^h(.),u(.)) | \le C_2 h \text{ uniformly for } u(.) \in \mathbf{U}_{ad,Lip}(L) \qquad \text{--- (3.12)}$$

for a constant $C_2$ to be calculated in the proof of this theorem.

Proof: We have

$$\begin{aligned}
&| J(x(.),u(.)) - J^h(x(.),u(.)) | = \\
&= \left| \int_0^T F(t,x(t),u(t))\,dt - \sum_{i=1}^{N-1} hF(ih,x(ih),u(ih)) \right| \le \\
&\le \sum_{i=1}^{N-1} \int_{ih}^{(i+1)h} \{ | F(t,x(t),u(t)) - F(ih,x(t),u(t)) | + \\
&+ | F(ih,x(t),u(t)) - F(ih,x(ih),u(ih)) | \}\, dt \le \\
&\le \{ M_F + L_F M_x \lambda_n + L_F L \} h \equiv C_3 h
\end{aligned} \qquad \text{--- (3.13)}$$

$$\begin{aligned}
&| J^h(x(.),u(.)) - J^h(x^h(.),u(.)) | \le \\
&\le \sum_{i=0}^{N-1} h\, | F(ih,x(ih),u(ih)) - F(ih,x^h(ih),u(ih)) | + \\
&+ | F_0(x(T)) - F_0(x^h(Nh)) | \le \\
&\le L_F C_1 h + L_{F_0} C_1 h \equiv C_4 h
\end{aligned} \qquad \text{--- (3.14)}$$

Thus

$$\begin{aligned}
&| J(x(.),u(.)) - J^h(x^h(.),u(.)) | \le | J(x(.),u(.)) - J^h(x(.),u(.)) | + \\
&+ | J^h(x(.),u(.)) - J^h(x^h(.),u(.)) | \le C_3 h + C_4 h \equiv C_2 h
\end{aligned} \qquad \text{--- (3.15)}$$

This proves the assertion of this theorem. ///

If $u^h(.)$ is a discrete-time control that satisfies

$|u^h((i+1)h) - u^h(ih)| \le Lh \ \forall i = 0,1,...,N-1$

we can construct a continuous-time control, which we shall denote by $\tilde{u}^h(.)$ by using linear interpolation, as follows:

$$\tilde{u}^h((i+\vartheta)h) := \vartheta\tilde{u}^h((i+1)h) + (1-\vartheta)\tilde{u}^h(ih)\ ,\ 0 \le \vartheta \le 1,\ 0 \le i \le N-1 \qquad \text{--- (3.16)}$$

We shall need the following:

<u>Lemma 3.1.</u> The function $\tilde{u}^h(.)$, costructed as in (3.16), is a member of $\mathbf{U}_{ad,Lip}(L)$.

<u>Proof:</u> We start with the following proposition: if $t_1 < t_2 < t_3$ and $u_1, u_2, u_3$ are such that $\left|\dfrac{u_2-u_1}{t_2-t_1}\right| \le L, \left|\dfrac{u_3-u_2}{t_3-t_2}\right| \le L$, if

$s_1 = \vartheta t_1 + (1-\vartheta)t_2,\ s_2 = \mu t_2 + (1-\mu)t_3, \vartheta \in [0,1],\ \mu \in [0,1],$

$\tilde{u}_1 = \vartheta u_1 + (1-\vartheta)u_2,\ \tilde{u}_2 = \mu u_2 + (1-\mu)u_3,$

then $\left|\dfrac{\tilde{u}_2-\tilde{u}_1}{s_2-s_1}\right| \le L$. This is proved by direct calculation:

$\tilde{u}_2 - \tilde{u}_1 = \mu u_2 + (1-\mu)u_3 - \vartheta u_1 - (1-\vartheta)u_2 = (1-\mu)(u_3-u_2) + \vartheta(u_2-u_1)$, and similarly $s_2 - s_1 = (1-\mu)(t_3-t_2) + \vartheta(t_2-t_1)$, thus

$|\tilde{u}_2 - \tilde{u}_1| \le (1-\mu)|u_3-u_2| + \vartheta|u_2-u_1| \le L[(1-\mu)(t_3-t_2) + \vartheta(t_2-t_1)] = L(s_2-s_1) =$
$= L|s_2-s_1|$

(the last equality is due to the fact that $s_2 - s_1 = (1-\mu)(t_3-t_2) + \vartheta(t_2-t_1)$ implies that $s_2 - s_1 \ge 0$).

This proposition can be extended, by induction, as follows:

If $t_1 < t_2 < ... < t_k < t_{k+1}$ and $u_1, u_2, ..., u_k, u_{k+1}$ are such that $\left|\dfrac{u_{j+1}-u_j}{t_{j+1}-t_j}\right| \le L\ \ \forall j = 1,2,...,k$, if

$s_1 = \vartheta_1 t_1 + (1-\vartheta_1)t_2,\ s_k = \vartheta_k t_k + (1-\vartheta_k)t_{k+1},$

$\tilde{u}_1 = \vartheta_1 u_1 + (1-\vartheta_1)u_2,\ \tilde{u}_k = \vartheta_k u_k + (1-\vartheta_k)u_{k+1},$

then $\left|\dfrac{\tilde{u}_k-\tilde{u}_1}{s_k-s_1}\right| \le L$. The inductive step is this: if the wanted conclusion is true for $k = m$, we shall show that it must be true for k=m+1; to that effect, we apply the proposition proved above, with the points $s_1, s_m, t_{m+1}$ in lieu of $t_1, t_2, t_3$, and the values $\tilde{u}_1, \tilde{u}_m, u_{m+1}$ instead of $u_1, u_2, u_3$, and we conclude that $\left| \frac{u_{m+1} - \tilde{u}_1}{t_{m+1} - s_1} \right| \le L$ ; using this last inequality, we apply again the same proposition to points $s_1, t_{m+1}, t_{m+2}$ in lieu of $t_1, t_2, t_3$, and values $\tilde{u}_1, u_{m+1}, u_{m+2}$ instead of $u_1, u_2, u_3$, and we obtain $\left| \frac{\tilde{u}_{m+1} - \tilde{u}_1}{s_{m+1} - s_1} \right| \le L$.

The assertion of the theorem is a particular case, with uniformly spaced points $t_j$, of the proposition that we just proved by induction for arbitrarily spaced points. ///

4. Approximate solution of the Volterra control problem with Lipschitz controls.

We consider the problem of minimizing the functional J given by (1.2) subject to the Volterra integral equation (1.1) and the constraint $u(.) \in \mathbf{U}_{ad,Lip}(L)$. We also consider the approximate problem of minimizing the fumctional $J^h$ given by (3.2) subject to the discretized Volterra equation (3.1) and the constraint $| u^h((i+1)h) - u^h(ih) | \le Lh \;\; \forall i = 0,1,...,N-2$. A solution of the discretized optimal control problem can be found by using the discrete dynamic programming equations of section 2.

We denote by $u^{h,*}(.)$ a solution of the discretized optimal control problem of this section. Our goal is to prove that $u^{h,*}(.)$ is close to an optimal control for the continuous optimal control problem of this section, in the following sense: if we construct a continuous-time control function by linear interpolation from the values of $u^{h,*}(.)$ and then use that continuous-time control function in the continous-time Volterra equation and the continuous-time functional J, then the value of J will be close to the infimum of J under the constraints stated above. Now, we make all this precise.

We denote by $\tilde{u}^{h,*}(.)$ the continuous-time control obtained through linear interpolation from $u^{h,*}(.)$ :

$$\tilde{u}^{h,*}(s) := \left(1 + i - \frac{s}{h}\right) u^{h,*}(ih) + \left(\frac{s}{h} - i\right) u^{h,*}((i+1)h) \;\text{ for } ih \le s \le (i+1)h,\; 0 \le i \le N-2 \qquad \text{--- (4.1)}$$

According to the results of section 3, $\tilde{u}^{h,*}(.)$ is in $\mathbf{U}_{ad,Lip}(L)$. We denote by $x^{h,*}(.)$ the solution of the discrete Volterra equation obtained by using control function $u^{h,*}(.)$, and by $x^*(.)$ the solution of the controlled Volterra integral equation (1.1) obtained by using the control function $\tilde{u}^{h,*}(.)$. We shall prove:

<u>Theorem 4.1.</u> We have

$$\lim_{h \to 0^+} J(x^*(.), \tilde{u}^{h,*}(.)) = \inf_{u(.) \in \mathbf{U}_{ad,Lip}(L)} J(x(.), u(.)) \qquad \text{--- (4.2)}$$

with linear rate of convergence, i.e.

$$\lim_{h\to 0^+} J(x^*(.),\tilde{u}^{h,*}(.)) - \inf_{u(.)\in \mathbf{U}_{ad,Lip}(L)} J(x(.),u(.)) = \mathbf{O}(h) \text{ as } h \to 0^+ .$$

<u>Proof:</u> We have, by the optimality of $u^{h,*}(.)$,

$$J^h(x^{h,*}(.),\tilde{u}^{h,*}(.)) \le J^h(x^h(.),u(.)) \ \ \forall u(.) \in \mathbf{U}_{ad,Lip}(L)$$

--- (4.3)

By combining (4.3) with theorem 3.2, we obtain

$$J^h(x^{h,*}(.),\tilde{u}^{h,*}(.)) \le J^h(x^h(.),u(.)) \le J(x(.),u(.)) + C_2h \ \forall u(.) \in \mathbf{U}_{ad,Lip}(L)$$

and

$$J(x^*,\tilde{u}^{h,*}(.)) \le J^h(x^{h,*}(.),\tilde{u}^{h,*}(.)) + C_2h \le J(x(.),u(.)) + 2C_2h \ \ \forall u(.) \in \mathbf{U}_{ad,Lip}(L)$$

thus, if we set $J^* := \inf_{u(.)\in \mathbf{U}_{ad,Lip}(L)} J(x(.),u(.))$ , we have

$$J^* \le J(x^*,\tilde{u}^{h,*}(.)) \le J^* + 2C_2h$$

--- (4.4)

which proves the assertion of the theorem. ///

5. Estimation of the computational cost of our version of dynamic programming, and design of parallel implementation.

The question of computational cost is important for every numerical algorithm; comparisons among different algorithms are generally based on their computational costs, although other aspects may also become relevant in specific cases.

We now present the calculation of the (approximate) computational cost of the variant of dynamic programming decribed in section 2.

First, at each step of numerically evaluating $V(i,\beta)$ by using (2.11), the variable β has to be quantized. If M is the number of quantized values of each $\beta(j)$, then the number of quantized values of β at the i-th stage of our dynamic programming is $M^i$. The number M depends on the dimension m of the range of the control and on the nature of the set **U**. If **U** is a cube, i.e. $\mathbf{U}=[a,b]^m$, and if Q is the quantized points in the interval $[a,b]$ in each coordinate of $\mathbb{R}^m$, then $M=Q^m$.
At the i-th stage of (2.11), and for each quantized value of β, the solution $x(i;i,\beta)$; we denote the cost of evaluating $x(i;i,\beta)$ by $\varphi_{sol}(i)$. In turn, $\varphi_{sol}(i)=i\,\varphi_{eval}(f;n)+\varphi_{eval}(x_0;n)+1$, where $\varphi_{eval}(f;n)$ is the cost of evaluating the values of the $\mathbb{R}^n$-valued function f, and likewise $\varphi_{eval}(x_0;n)$ is the cost of evaluating $x_0(i)$; in order to evaluate $x(i;i,\beta)$ by using (2.1), i evaluations of f, one evaluation of $x_0$, and one multiplication by h, are required. The form of $\varphi_{sol}(i)$ reflects the fact that the total number of evaluations of $x_0$, over the range $0\le i\le N$, is N+1. (As usual for work of this type, we do not include the operations of addition and subtraction in the estimation of the computational cost.) The optimization shown on the right-hand side of (2.11) has a computational cost which we denote by $\varphi_{opt}(i)$. In order to produce a result that can be useful for actual calculations, a minimizer $\xi^*(i)$ of the expression on the right-hand side of (2.11) would have to be approximately expressed, for instance via linear interpolation, as a function of quantized values of the control; we denote this cost of interpolation by $\varphi_{int}(\xi^*(i))$. Thus far we have expressed these costs in connection with the evaluation of $V(i,\beta)$. We denote by $\varphi_{eval}(V;i)$ the cost of evaluating $V(i,\beta)$ over all quantized values of β, and by $\varphi_{eval}(\Phi;\varphi_{sol}(i),\xi)$ the cost of evaluating Φ for each β and ξ, then

$$\varphi_{eval}(V;i)=\varphi_{eval}(V;i+1)+M^{i+1}\varphi_{eval}(\Phi;\varphi_{sol}(i),\xi)+$$
$$+M^i[\varphi_{opt}(i)+\varphi_{int}(\xi^*(i))]$$

--- (5.1)

Eq. (5.1) is supplemented by the final-time condition

$$\varphi_{eval}(V;N) = M^{N+1}\varphi_{eval}(\Phi;\varphi_{sol}(N),\xi)$$

--- (5.2)

In order to arrive at closed-form estimates, we make certain simplifying assumptions. First, we assume that $\varphi_{opt}(i)$ and $\varphi_{int}(\xi^*(i))$ are constants, and we set $A := \varphi_{opt}(i) + \varphi_{int}(\xi^*(i))$; second, we assume that the cost of evaluating $\Phi$ is proportional to the cost of evaluating x(i;i,β), thus we assume that $\varphi_{eval}(\Phi;\varphi_{sol}(i),\xi) = C_{eval}[i\varphi_{eval}(f;n) + \varphi_{eval}(x_0;n) + 1] \equiv C_{\Phi,1}i + C_{\Phi,0}$.

Then the system consisting of (5.1) and (5.2) becomes

$$\varphi_{eval}(V;i) = \varphi_{eval}(V;i+1) + M^{i+1}(C_0 + iC_1) + M^i A;$$
$$\varphi_{eval}(V;N) = M^{N+1}(C_0 + NC_1)$$

--- (5.3)

The total cost of evaluating V is

$$\varphi_{total}(V) = \sum_{i=1}^{N} \varphi_{eval}(V;i)$$

--- (5.4)

The summations indicated in (5.3) and (5.4) can be carried out by elementary methods, and the answer is

$$\varphi_{total}(V) = C_0 \frac{NM^{N+2} - (N+1)M^{N+1} + M}{(M-1)^2} +$$
$$+ C_1 \frac{N^2M^{N+3} - (2N^2 + 2N - 1)M^{N+2} + (N+1)^2 M^{N+1} - M^2 - M}{(M-1)^3} +$$
$$+ A_0 \frac{(N-1)M^{N+1} - NM^N + M}{(M-1)^2}$$

--- (5.5)

The important conclusion follows from the highest-degree term in (5.5): under the stated assumptions, the computational cost of our dynamic programming algorithm is of the order of $N^2M^N$; for a rectangular set **U**, the cost is of the order of $N^2Q^{mN}$.

Next, we examine the possibilities of parallel implementation of the discrete dynamic programming equations for Volterra control problems. At each stage (discrete time) i, the evaluation of $x(i;i,\beta)$ and the evaluation of the right-hand side of the discrete dynamic programming equation (2.25) can be carried out in parallel for each of the quantized values of β. Thus, this type of calculation is suited to SIMD (single instruction, multiple data) parallel processors. At each stage, and for each quantized value of β, all values of $V(i+1,\beta\otimes\xi)$ over all quantized values of ξ, are needed. These values, over all quantized values of β and ξ, have to be stored in a shared memory unit with which all processors can communicate and select those values of β that correspond to the appropriate processor. The number of processors that are needed at each stage i is a function of i, since the dimensionality of β depends on i. This leads to an adaptive requirement: the number of active processors is time-varying, it depends on the discrete time i at which the set of parallel computations needs to be performed. This adaptivity is a normal feature of parallel compuing, cf. [A]. The cost, in terms of computing time rather than number of operations, at each stage i, is $\varphi_{eval}(V;i)+\varphi_{comm}(i)+\varphi_{sel}(i)$, where $\varphi_{comm}(i)$ is the cost of the processors' communicating with the memory unit at stage i, and $\varphi_{sel}(i)$ the cost associated with selecting the values of β, out of all $V(i+1,\beta\otimes\xi)$, that correspond to each processor. In general, the number of processors will be smaller that the number of quantized values of the control function, thus a set of values of the control would need to be assigned to each processor.

## 6. Remarks on continuous-time dynamic programming for Volterra control.

This section concerns a conceptual question that arises naturally from the results of the previous section, namely: what, if any, would be the form of continuous-time dynamic programming equations for optimal control of systems governed by Volterra integral equations? This question is, as far as we can judge, of only conceptual value: any actual computational solution of Volterra control problems will require some sort of approximation, such as the method we have developed above. It is nevertheless a question that a reader might reasonably ask.

It is expressly stated that this section does not contain rigorous results; in our assessment, a rigorous development would be useless for actually solving the related optimal control problems. Our purpose here is to formally discern the possible nature of continuous-time dynamic programming equations for systems governed by Volterra integral equations, not to prove theorems about such equations.

The continuous-time analogue of (2.11) is not a straightforward matter. In the continuous case, if we use, as a variable in the Hamilton-Jacobi function, the restriction $u(.)|_{[0,t]}$ of a control function to the interval [0, t], then the pairs $(t, u(.)|_{[0,t]})$ form a vector bundle, each $u(.)|_{[0,t]}$ being an element of the space $C([0,t] \mapsto \mathbb{R})$, the space of continuous functions from [0, t] into the real numbers. Because the spaces $C([0,t] \mapsto \mathbb{R})$ depend on t, differentiation with respect to the variables t and $u(.)|_{[0,t]}$, in the Hamilton-Jacobi function, cannot be carried separately with respect to each variable, and thus there is no straightforward way to obtain, even formally, a continuous-time dynamic programming equation.

We have devised a roundabout way to obtain a framework that allows differentiation, by using a transformation that changes the variable vector spaces into a fixed vector space, and we carry out the calculations in these transformed spaces.
For each function v(.) in $C([0,t] \mapsto \mathbb{R})$, we define the corresponding function $v_t^{\#}(.)$ in $C([0,1] \mapsto \mathbb{R})$ by $v_t^{\#}(\tau) = v(t\tau)$, $0 \le \tau \le 1$. Similarly, for each continuous function of two variables, say $w(t_1, s_1)$, defined for $0 \le s_1 \le t_1 \le t$, we define the corresponding function $w_t^{\#}(.,.)$ in $C(\Delta([0,1]) \mapsto \mathbb{R})$, where $\Delta([0,1]) := \{(\tau,\sigma) : 0 \le \sigma \le \tau \le 1\}$, by $w_t^{\#}(\tau,\sigma) := w(t\tau, t\sigma)$.
Now, the integral equation

$$x(t) = x_0(t) + \int_0^t f(t, s, x(s), u(s))ds$$

can be written in the form

$$x(t\tau) = x_0(t\tau) + \int_0^{\tau} f(t\tau, t\sigma, x(t\sigma), u(t\sigma))\, t\, d\sigma$$

--- (6.1)

which is the same as

$$x_t^{\#}(\tau) = x_{0,t}^{\#}(\tau) + t\int_0^{\tau} f_t^{\#}(\tau, \sigma, x_t^{\#}(\sigma), u_t^{\#}(\sigma))\, d\sigma,\ 0 \le \tau \le 1$$

--- (6.2)

The cost functional

$$J = \int_0^{T} F(t, x(t), u(t))\, dt + F_0(x(T))$$

can also be written as

$$J = \Psi(x_T^{\#}, u_T^{\#}) := \int_0^{1} F(T\tau, x_T^{\#}(\tau), u_T^{\#}(\tau))\, T\, d\tau + F_0(x_T^{\#}(1))$$

--- (6.3)

We denote by **R** the solution operator associated with (6.2), i.e.

$$x_t^{\#} = \mathbf{R}(t, u_t^{\#})$$

--- (6.4)

We define the Hamilton-Jacobi function $V(t, \alpha^{\#})$ by

$$V(t, \alpha^{\#}) := \inf\{\Psi(x_T^{\#}, u_T^{\#}) : u_t^{\#} = \alpha^{\#}\}$$

--- (7.5)

We consider controls $u_t^{\#}$ that satisfy $\frac{\partial u_t^{\#}}{\partial t} = w_t^{\#}$ for some continuous function $w_t^{\#}$, and we interpret $w_t^{\#}$ as a new control taking values in a set W of continuous functions. Then the dynamic programming equation is, formally,

$$\frac{\partial V(t,u^{\#})}{\partial t} + \inf_{w^{\#} \in W} \left\langle \frac{\partial V(t,u^{\#})}{\partial u^{\#}}, w^{\#} \right\rangle = 0; \; V(T,u^{\#}) = \Psi(R(T,u^{\#}),u^{\#}) \qquad \text{--- (6.6)}$$

The brackets $\langle \cdot, \cdot \rangle$ are used to signify the action of the linear functional $\frac{\partial V(t,u^{\#})}{\partial u^{\#}}$ on $w^{\#}$.